\def\draft{n}
\documentclass{amsart}
\usepackage[headings]{fullpage}
\usepackage{amssymb, epsfig,amsbsy,amsmath,bbold}

\theoremstyle{plain}

\newtheorem{theorem}{Theorem}
\newtheorem{proposition}{Proposition}[section]

\theoremstyle{definition}

\theoremstyle{remark}

\newtheorem{remark}[proposition]{Remark}

\def\printname#1{
        \if\draft y
                \smash{\makebox[0pt]{\hspace{-0.5in}
                        \raisebox{8pt}{\tt\tiny #1}}}
        \fi
}

\newcommand{\psdraw}[2]
         {\begin{array}{c} \hspace{-1.3mm}
        \raisebox{-4pt}{\epsfig{figure=draws/#1.eps,width=#2}}
        \hspace{-1.9mm}\end{array}}

\newlength{\standardunitlength}
\catcode`\@=11
\long\def\@makecaption#1#2{%
     \vskip 10pt

\setbox\@tempboxa\hbox{
       \small\sf{\bfcaptionfont #1. }\ignorespaces #2}%
     \ifdim \wd\@tempboxa >\captionwidth {%
         \rightskip=\@captionmargin\leftskip=\@captionmargin
         \unhbox\@tempboxa\par}%
       \else
         \hbox to\hsize{\hfil\box\@tempboxa\hfil}%
     \fi}
\font\bfcaptionfont=cmssbx10 scaled \magstephalf
\newdimen\@captionmargin\@captionmargin=2\parindent
\newdimen\captionwidth\captionwidth=\hsize
\catcode`\@=12

\def\lbl#1{\label{#1}\printname{#1}}

\def\BN{\mathbb N}
\def\BZ{\mathbb Z}

\def\BQ{\mathbb Q}

\def\BC{\mathbb C}

\def\la{\langle}
\def\ra{\rangle}

\def\SL{\mathrm{SL}}

\def\Cert{\mathrm{Cert}}
\def\hatJ{\hat{J}}

\newcommand{\qfac}[1] {(q; q)_{#1}}
\newcommand{\qrf}[2] {({#1}; q)_{#2}}

\def\checkJ{\check{J}}

\begin{document}


\title[The non-commutative $A$-polynomial of twist knots]{
The non-commutative $A$-polynomial of twist knots}

\author{Stavros Garoufalidis}
\address{School of Mathematics \\
         Georgia Institute of Technology \\
         Atlanta, GA 30332-0160, USA \\ 
         {\tt http://www.math.gatech} \newline {\tt .edu/$\sim$stavros } }
\email{stavros@math.gatech.edu}
\author{Xinyu Sun}
\address{Department of Mathematics \\
        Tulane University \\
        6823 St. Charles Ave \\
        New Orleans, LA 70118, USA \\
{\tt http://www.} \newline {\tt math.tulane.edu/$\sim$xsun1}}
\email{xsun1@tulane.edu}

\thanks{S.G. was supported in part by National Science 
Foundation. \\
\newline
1991 {\em Mathematics Classification.} Primary 57N10. Secondary 57M25.
\newline
{\em Key words and phrases: Knots, Jones polynomial, colored Jones function,
$A$-polynomial, $C$-polynomial, non-commutative $A$-polynomial, 
$q$-difference equations, WZ Algorithm, Creative Telescoping, 
Gosper's Algorithm, certificate, multi-certificate.
}
}

\date{July 6, 2009.}


\begin{abstract}
The purpose of the paper is two-fold: to introduce a multivariable creative
telescoping method, and to apply it in a problem of Quantum Topology: namely
the computation of the non-commutative $A$-polynomial of twist knots.

Our multivariable creative telescoping
method allows us to compute linear recursions for sums of the form
$J(n)=\sum_k c(n,k) \hatJ (k)$ given a recursion relation for $(\hatJ(n))$
a the hypergeometric kernel $c(n,k)$. As an application of our method,
we explicitly compute the
non-commutative $A$-polynomial for twist knots with $-8$ and $11$ crossings.

The non-commutative $A$-polynomial of a knot encodes the monic, linear,
minimal order $q$-difference equation satisfied by the sequence of
colored Jones polynomials of the knot. 
Its specialization to $q=1$ is conjectured to
be the better-known $A$-polynomial of a knot, which encodes important
information about the geometry and topology of the knot complement.
Unlike the case of the Jones polynomial, which is easily computable for
knots with $50$ crossings, the $A$-polynomial is harder to compute and 
already unknown for some knots with $12$ crossings. 
\end{abstract}

\maketitle

\tableofcontents


\section{Introduction}
\lbl{sec.intro}

\subsection{The goal}
\lbl{sub.goal}

The purpose of the paper is two-fold: to introduce a multivariable creative
telescoping method, and to apply it in a problem of Quantum Topology: namely
the computation of the non-commutative $A$-polynomial of twist knots.

Our multivariable creative telescoping
method allows us to compute linear recursions for sums of the form
$J(n)=\sum_k c(n,k) \hatJ (k)$ given a recursion relation for $(\hatJ(n))$
a the hypergeometric kernel $c(n,k)$. General theory implies the existence
of a recursion relation for $(\hatJ(n))$. However, in practice the computation
is not manageable for twist knots, and there is no guarantee that the recursion
relation will be of minimal order. 
Our method does not guarantee a minimal order recursion 
relation either, however (unlike the known methods) it is manageable 
and produces a minimal order recursion relation the non-commutative 
$A$-polynomial for twist knots of $-8, \dots, 11$ twists. 
The non-commutative $A$-polynomial encodes the unique 
monic, linear, minimal order $q$-difference equation satisfied 
by the sequence of colored Jones polynomials of the knot. 
Our results give a new proof of the AJ-Conjecture for those knots.

\subsection{The Jones polynomial of a knot}
\lbl{sub.history}

In this section we recall the relevant Laurent polynomial invariants of
knots, such as the Jones polynomial and its colored cousins. 
In 1985 V. Jones introduced the famous {\em Jones polynomial} of a knot
$K$ in 3-space, \cite{Jo}. The Jones polynomial (an element of 
$\BZ[q^{\pm 1}]$)
is a powerful knot invariant which amongth other things detects cheirality,
and it can be extended to a sequence $(J_K(n))$
of Laurent polynomials by taking parallels of a knot $K$. Technically,
$J_K(n) \in \BZ[q^{\pm 1}]$ is the quantum group invariant of the 0-framed knot
using the $n$-dimensional representation of $\mathrm{SU}(2)$, and normalized
by $J_{\text{Unknot}}(n)=1$. For a detailed definition, see \cite{Tu} and also
\cite{GL1}. With this normalization, we have that $J_K(1)=1$,
and $J_K(2)$ is the Jones polynomial of $K$. 

For a given knot $K$, the sequence of Laurent polynomials $(J_K(n))$
is not random. To be precise, $(J_K(n))$ is $q$-{\em holonomic} i.e., it
satisfies a linear $q$-difference equation (which of course, depends on 
the knot) with coefficients in $\BQ(q,q^n)$. 
This fact, proven in \cite{GL1}, is an easy consequence of two facts:

\begin{itemize}
\item[(a)] 
$J_K(n)$ is a finite multisum of a proper $q$-hypergeometric term,
as follows from the state-sum definition of the colored Jones function;
see \cite{GL1}.
\item[(b)]
ultisums of proper $q$-hypergeometric terms are $q$-holonomic, as follows
from the WZ theory of Wilf-Zeilberger; see \cite{WZ}.
\end{itemize}

\subsection{The non-commutative $A$-polynomial of a knot and its significance}
\lbl{sub.Aq}

A $q$-holonomic sequence is annihilated by a unique monic homogeneous 
linear $q$-difference equation of smallest degree, and the corresponding 
monic polynomial in two $q$-commuting variables $E$ and $Q$ is an invariant 
(the so-called {\em characteristic polynomial} of the $q$-holonomic
sequence. We define the {\em non-commutative $A$-polynomial} 
$A_K(E,Q,q)$ of a knot $K$ to be the characteristic polynomial of 
$(J_K(n))$. 

In \cite{Ga}, it was conjectured by the first author 
(the so-called AJ Conjecture) that the specialization $A_K(E,Q,1)$ of
$A_K(E,Q,q)$ should agree with the {\em $A$-polynomial} of a knot. 
The latter is
an important invariant that parametrizes the $\SL(2,\BC)$ character variety
of the knot complement, as viewed from the boundary torus. For a detailed
definition of the $A$-polynomial, its properties and its applications
to the geometry and topology of the knot complement, see \cite{CCGLS}.
Thus, $A_K(E,Q,q)$ can be thought of as a deformation (or quantization) of
the character variety.  

The Jones polynomial of a knot is easily computable via skein theory with
knots with, say, $50$ crossings; see for example \cite{B-N}. On the other hand,
the $A$-polynomial of a knot is much harder to compute, and at present it
is unknown for some knots with $12$ crossings. There are two general
methods to compute the $A$-polynomial: an exact (primarily elimination,
and Puiseux expansions) developed by Boyd \cite{Bo} and a numerical one
developed by Culler \cite{Cu}.

The non-commutative $A$-polynomial and its possible relation with the 
$A$-polynomial of a knot is an important ingredient to the
Hyperbolic Volume Conjecture and its generalization.

\subsection{Computing the non-commutative $A$-polynomial}
\lbl{sub.comuteAq}

For theoretical as well as experimental reasons it would be good to have 
explicit formulas for the non-commutative $A$-polynomial. 
So far, an explicit formula has been
given for torus knots in \cite{Ge} (using properties of the Kauffman 
bracket skein module of the solid torus), as well as for the simplest 
hyperbolic $4_1$ knot in \cite{GL1} (using an explicit single-sum formula
for the colored Jones function).

The WZ algorithm has been implemented (see \cite{PWZ,PR2,PR3}) and together 
with
explicit state-sum formulas for the colored Jones function of an arbitrary
planar projection given in \cite{GL1}, in principle one can obtain a 
linear $q$-difference equation for the colored Jones function of an arbitrary
knot. There are two problems with this approach:

\begin{itemize}
\item[(a)] 
The number of summation variables in the multisum formulas is
generally two less than the number of crossings, and the $q$-multisum
algorithms appear to be slow for the current machines.
\item[(b)]
There is no guarantee that the various $q$-multisum algorithms will give
a minimal order linear $q$-difference equation. In fact, in many cases
(where symmetry is involved), it has been observed that they fail to
give the minimal order $q$-difference equation. See \cite{PR1} for 
well-known examples of this failure.
\end{itemize}
With respect to the first problem, we were unable to use the sofware of
\cite{AZ,S} to compute a $q$-difference equation for our double sums.
  
One might wonder whether Problem (b) really occurs for the state-sums
that originate in knot theory. As expected, this problem {\em does} occur. 
The knots $5_2$ and $6_1$ have double-sum formulas for their colored Jones 
function. An application of the $q$-multisum package of \cite{PR3} was done by
Takata in \cite{Ta} who found out an explicit inhomogeneous $q$-difference 
equation of degree $5$ and $5$ respectively. On the other hand, as we shall 
see, there exist inhomogeneous $q$-difference equations of degree $3$
and $4$ respectively.

In a different direction, Le used geometric methods of the Kauffman
bracket skein module and was able to prove the AJ Conjecture for most 
2-bridge knots, as well as give a linear algebra algorithm that in principle 
computes the non-commutative $A$-polynomial; see \cite[Thm.1]{Le} and 
\cite[Sec.5.6.3]{Le}. The algorithm was implemented
in Maple by the second author, but proved to be too slow to run
for the $5_2$ and $6_1$ knots.

\subsection{A sample of our results}
\lbl{sub.results}

The main goal of the paper is to give an explicit formula for the 
non-commutative $A$-polynomial of twist knots with $p$ twists, where
$p=-8,\dots,11$. 

Let us recall the {\em twist knots} $K_p$ for integer $p$, shown 
in Figure \ref{twist}. The planar projection of $K_p$ 
has $2|p|+2$ crossings, $2|p|$ of which come from $p$ full twists, 
and $2$ come from the negative {\em clasp}.

\begin{figure}[htpb]
$$ 
\psdraw{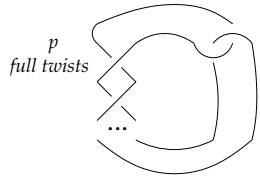}{1.3in} 
$$
\caption{The twist knot $K_p$, for integers $p$.}\lbl{twist}
\end{figure}

For small $p$, these knots may be identified with ones from Rolfsen's
table (see \cite{Rf} and \cite{B-N}) as follows:
$$
K_1=3_1, \quad K_2=5_2, \quad K_3=7_2, \quad K_4=9_2,
$$
$$
K_{-1}=4_1, \quad K_{-2}=6_1, \quad K_{-3}=8_1, \quad K_{-4}=10_1.
$$
Let $E$ and $Q$ denote the operators that act on a sequence $(J(n))$
of Laurent polynomials $J(n) \in \BZ[q^{\pm 1}]$ by:

\begin{equation}
\lbl{eq.EQ}
(EJ)(n)=J(n+1), \qquad (QJ)(n)=q^n J(n).
\end{equation}
Note that $EQ=q Q E$. 
Let $(A^{nh}_p(E,Q,q),B_p(q^n,q))$ denote the {\em inhomogeneous} 
non-commutative $A$-polynomial of $K_p$. That is, $A^{nh}_p(E,Q,q)$
is monic and minimal degree (with respect to $E$) that satisfies the 
equation

\begin{equation}
\lbl{eq.Anh}
A^{nh}_p(E,Q,q) J_p=B_p(q^n,q),
\end{equation}
where $B_p(q^n,q) \in \BQ(q^n,q)$. To convert the inhomogeneous equation
above to a homogeneous one, see Section \ref{sec.reviewct}.


\begin{theorem}
\lbl{thm.0}
\rm{(a)}
For $p=\pm 1$ we have: 
{\small  
\begin{eqnarray*}
A^{nh}_1(E,Q,q)	
&=& 
q^{3n+2}(q^n-1)+(q^{n+1}-1)E
\\
B^{nh}_1(q^n,q)	
&=&
(q^{2n+1}-1)q^n
\\
A^{nh}_{-1}(E,Q,q)	
&=& 
q^{2n+2}(q^n-1)(q^{2n+3}-1)
-(q^{n+1}+1)(q^{4n+4}-q^{3n+3}-q^{2n+3}-q^{n+1}-q^{2n+1}+1)(q^{n+1}-1)^2E
\\
& & +q^{2n+2}(q^{2n+1}-1)(q^{n+2}-1)E^2
\\
B^{nh}_{-1}(q^n,q)	
&=&
q^{n+1}(q^{2n+3}-1)(q^{n+1}+1)(q^{2n+1}-1)
\end{eqnarray*}
}
\rm{(b)} For $p=\pm 2$ we have: 
{\small
\begin{eqnarray*}
\lbl{lbl.eq.2}
A^{nh}_2(E,Q,q) &=& 
q^{7n+9}(-1+q^n)(q^{2n+4}-1)(q^{2n+5}-1)-q^{2n+5}(q^{n+1}-1)(q^{2n+2}-1)(q^{2n+5}-1)\\
\nonumber
& & (q^{5n+6}-2q^{4n+5}-q^{3n+5}+q^{2n+4}-q^{2n+3}-q^{2n+2}-q^{3n+2}+q^{2n+1}+q^{n+1}-1)E\\
\nonumber
& & +q(q^{n+2}-1)(q^{2n+1}-1)(q^{2n+4}-1)(q^{5n+9}-q^{4n+7}-q^{3n+7}+q^{3n+6}+q^{3n+5}\\
\nonumber
& & +q^{2n+5}-q^{3n+4}+2q^{n+2}+q^{2n+2}-1)E^2+(q^{n+3}-1)(q^{2n+1}-1)(q^{2n+2}-1)E^3,\\
\nonumber
B_2(q^n,q) &=& (q^{n+1}+1)(q^{2n+1}-1)(q^{n+2}+1)(q^{2n+5}-1)(q^{2n+3}-1)q^{2n+4},\\	
\nonumber\lbl{lbl.eq.-2}
A^{nh}_{-2}(E,Q,q)	&=& 
(-1+q^n)(q^{2n+5}-1)(q^{2n+6}-1)(q^{2n+7}-1)q^{4n+8}-q^{2n+5}(q^{n+1}-1)(q^{2n+6}-1)\\
& & (q^{2n+7}-1)(q^{2n+2}-1)(q^{6n+9}-q^{5n+8}-q^{4n+8}+q^{4n+7}\\	
\nonumber
& & +q^{3n+7}+q^{4n+6}-q^{3n+6}-q^{2n+6}-q^{3n+5}-q^{2n+5}-q^{4n+4}+q^{3n+3}\\	
\nonumber
& & +q^{2n+3}-q^{n+3}-q^{2n+2}-2q^{n+2}-q^{2n+1}+q+1)E\\	
\nonumber
& & +q(q^{n+2}-1)(q^{2n+1}-1)(q^{2n+4}-1)(q^{2n+7}-1)
(1-q^{2n+3}-2q^{n+2}-q^{n+3}+q^{2n+5}\\ 
\nonumber
& & +2q^{5n+8}
+3q^{4n+8}-q^{3n+7}-q^{4n+6}-2q^{3n+6}-q^{2n+6}+q^{3n+5}-q^{2n+2}+q^{4n+5}+q^{2n+4}\\	
\nonumber
& & -q^{6n+10}-q^{6n+14}-q^{5n+11}-q^{7n+15}-q^{6n+11}+q^{6n+12}+q^{8n+16}+q^{5n+13}\\
\nonumber
& & +q^{6n+13}+q^{5n+9}+q^{3n+9}+2q^{3n+4}-2q^{7n+14}+3q^{4n+9}+2q^{3n+8}-2q^{5n+10}\\
\nonumber
& & +2q^{5n+12})E^2-(q^{n+3}-1)(q^{2n+6}-1)(q^{2n+2}-1)(q^{2n+1}-1)(q^{6n+16}+q^{6n+15}\\
\nonumber
& & -q^{5n+14}-2q^{5n+13}-q^{4n+13}-q^{4n+12}+q^{4n+10}+q^{3n+10}-q^{4n+9}-q^{3n+9}\\
\nonumber
& & -q^{4n+8}-q^{3n+8}-q^{2n+7}+q^{3n+6}+q^{2n+6}+q^{2n+5}-q^{2n+3}-q^{n+3}+1)E^3\\
\nonumber
& & +q^{4n+12}(q^{n+4}-1)(q^{2n+1}-1)(q^{2n+2}-1)(q^{2n+3}-1)E^4,\\
\nonumber
B_{-2}(q^n,q) &=& 
(q^{n+3}+1)(q^{2n+5}-1)(q^{2n+7}-1)q^{2n+6}(q^{n+2}+1)(q^{n+1}+1)(q^{2n+1}-1)(q^{2n+3}-1).
\end{eqnarray*}
}
\end{theorem}
The formulas quickly become too lengthy to type. For more information,
see Appendix \ref{sec.p-33} for $p=\pm 3$ as well as the data file
\cite{GS2} for $p=-8,\dots,11$.  
Theorem \ref{thm.0} gives a new proof of the AJ-Conjecture 
for twist knots with $-8,\dots,11$ twists.

\subsection{Plan of the proof}
In Section \ref{sec.strategy} we outline the main strategy. 
The idea is to use the recursion relation of the cyclotomic function 
$(\hatJ_p(n))$ of the twist knot $K_p$ (from \cite{GS1}) as well as
a single-sum relation between $J_p(n)$ and $\hatJ_p(n)$, together with
some new ideas of Creative Telescoping and some guessing.
In Section \ref{sec.reviewct} we review the method of Creative Telescoping
and in Section \ref{sec.hybrid} we present a multi-certificate version that 
takes
into account the product of a hypergeometric summand with a $q$-holonomic
one.

We conclude with three appendices: in Appendix
\ref{sec.GF} we present an alternative method that uses generating
functions (that was kindly communicated to us by Zeilberger). In 
Appendix \ref{sec.p-33} we give the non-commutative polynomial of twist knots 
$K_p$ for $p=-3,3$, and in Appendix \ref{sec.Acomm} we give the $A$-polynomial
of the same knots.

\subsection{Acknowledgement}
The authors wish to thank D. Zeilberger for encouragement and enlightening
conversations.

\section{The strategy}
\lbl{sec.strategy}

In this section we will describe our strategy to obtain a formula for 
$A_p(E,Q,q)$.

\begin{itemize}
\item[(a)]
We consider the cyclotomic function $\hatJ_K(n) \in \BZ[q^{\pm 1}]$
introduced by Habiro in \cite{Ha}, who used the notation $J_K(P''_n)$. 
\item[(b)]
The relation between the cyclotomic and the colored Jones functions
is given by:
\begin{equation}
\lbl{eq.CJ}
J_K(n) = \sum^n_{k=0} c(n, k) \hatJ_K(k),
\end{equation}
where the {\em cyclotomic kernel} $c(n,k)$ is a proper $q$-hypergeometric
term given for $0 \leq k \leq n$ by: 
\begin{eqnarray}
\lbl{eq.cnk}
c(n, k) &=& \frac{\{n-k\}\{n-k+1\} \dots \{n+k-1\}\{n+k\}}{\{n\}} 
\\ \notag
&=&
(-1)^k q^{-k(k+1)/2} \qrf{q^{1-n}}{k} \qrf{q^{1+n}}{k},
\end{eqnarray}
where 
$$
\{n\}=q^{n/2}-q^{-n/2},
$$
and the {\em quantum factorial} is defined by:
$$
(x; q)_n = \left\{
\begin{tabular}{ll}
$(1-x) \cdots (1-xq^{n-1})$ &   if $n>0$;  \\[1mm]
1       &       if $n = 0$;     \\
$\frac{1}{(1-xq^{-1}) \cdots (1-xq^{n})}$ & if $n < 0$.
\end{tabular}
\right. 
$$
Please note that we are using the unbalanced quantum factorials 
(common in discrete math) and not the balanced ones 
(common in the representation theory of quantum groups).
\item[(c)] 
The cyclotomic function $(\hatJ_K(n))$ is $q$-holonomic, 
as shown in \cite{GL1}, and its characteristic polynomial $C_K(E,Q,q)$
is defined to be the non-commutative $C$-polynomial of a knot.
Let us abbreviate $\hatJ_{K_p}(n)$, $J_{K_p}(n)$, $A_{K_p}(E,Q,q)$
and $C_{K_p}(E,Q,q)$
for the twist knots $K_p$ by $\hatJ_{p}(n)$, $J_{p}(n)$,
$A_{p}(E,Q,q)$ and $A_{p}(E,Q,q)$ respectively.
\item[(d)] 
In \cite{GS1} we gave an explicit formula for $C_p(E,Q,q)$.
\item[(e)]
Using the explicit formula for $C_p(E,Q,q)$ as well as the relation 
\eqref{eq.CJ} and a version of creative telescoping (and some guessing),
we deduce a linear $q$-difference equation for $(J_p(n))$, which specializes
to the $A$-polynomial of $K_p$ when $q=1$.
\item[(f)]
Since the $A$-polynomial of $K_p$ is irreducible (see \cite{HS1}),
and the non-commutative $A$-polynomial of $K_p$ specializes to the 
$A$-polynomial, 
it follows that our $q$-difference equation is indeed of minimal order.
This computes $A_p(E,Q,q)$.
\end{itemize}

For completeness and concreteness, we give a formula for $\hatJ_p(n)$, using 
\cite[Thm.5.1]{Ma} (compare also with \cite[Sec.3]{Ga}):

\begin{eqnarray}
\lbl{eq.Cpn}
\hatJ_p(n)&=& \sum_{k=0}^{n}     
        q^{n(n + 3)/2 + p k(k + 1) + k(k - 1)/2} (-1)^{n + k + 1}
        \frac{(q^{2k + 1} - 1)\qfac{n}} {\qfac{n + k + 1} \qfac{n-k}} 
\\
\notag
&=&
\sum_{k=0}^{n}     
        q^{n(n + 3)/2 + p k(k + 1) + k(k - 1)/2} (-1)^{n + k + 1}
        \frac{(q^{2k + 1} - 1)(q^{n-k+1};q)_k} {\qfac{n + k + 1}}.
\end{eqnarray}
Observe that since $(q^{n-k+1};q)_k=0$ for $k>n>0$, we can assume that
the $k$-summation in the above equation is for $0 \leq k < +\infty$.

Equations \eqref{eq.CJ}, \eqref{eq.cnk} and \eqref{eq.Cpn} imply that
$J_p(n)$ is given by a double-sum formula of a proper $q$-hypergeometric
summand. As explained earlier, the {\tt qMultisum.m} implementation of
the WZ algorithm given in \cite{PR2,PR3} and used in \cite{Ta}
is slow to run, and gives $q$-difference equations of higher than actual 
degree. An application of our multicertificate version of Creative
Telescoping is the following theorem.

\begin{theorem}
\lbl{thm.1}
The minimal inhomogeneous recursion for $J_p(n)$ for 
$-8 \le p \le 11$ is an explicit linear $q$-difference equation
of order
$$
\begin{cases}
2p-1 & \text{if} \qquad 0 < p \leq 11; \\
2|p|  & \text{if} \qquad -8 \leq p < 0.
\end{cases}
$$
\end{theorem}
The inhomogeneous recursion is given by Theorem \ref{thm.0} for 
$p=-2,\dots,2$, Appendix \ref{sec.p-33} for $p=\pm 3$
and the data file \cite{GS2} for $p=-8, \dots, 11$.

\section{A brief review of creative telescoping}
\lbl{sec.reviewct}

In this section we recall briefly some key ideas of Zeilberger on recursion
relations of combinatorial sums. An excellent reference is \cite{PWZ}.
For a longer introduction, see also \cite[Sec.3]{GS1}.

A {\em term} is $F(n, k)$ called {\em hypergeometric} 
if both $\frac{F(n+1, k)}{F(n, k)}$ and $\frac{F(n, k+1)}{F(n, k)}$
are rational functions over $n$ and $k$. In other words,

\begin{equation}
\lbl{eq.rationalF}
\frac{F(n+1, k)}{F(n, k)} \in \BQ(n,k), \qquad \frac{F(n, k+1)}{F(n, k)}
\in \BQ(n,k).
\end{equation}
Examples of hypergeometric terms
are $F(n,k)=(an+bk+c)!$ (for integers $a,b,c$), and ratios of products of such.
The latter are actually called {\em proper hypergeometric}.
A key problem is to construct recursion relations for sums of the form:

\begin{equation}
\lbl{eq.sumS}
S(n)=\sum_k F(n,k),
\end{equation}
where $F(n,k)$ is a proper hypergeometric term.
The summation set can be the set of all integers or an interval thereof.
Let us first suppose that summation is over entire set of integers.
Sister Celine \cite{Fa} (see also \cite{PWZ}) proved the following:

\begin{theorem}
\lbl{thm.celine}
Given a proper hypergeometric term $F(n, k)$, 
there exist a natural number $I \in \BN$ 
and a set of functions $a_i(n) \in \BQ(n)$, 
$0 \le i \le I$, such that

\begin{eqnarray} 
\lbl{eqn.celine}
\sum_{i=0}^{I} a_i(n) F(n+i, k) = 0.
\end{eqnarray}
\end{theorem}
The important part of the above theorem is that the functions $a_i(n)$ are 
independent of $k$. Therefore if we take the sum over $k$ on both sides, 
we get

\begin{eqnarray} 
\lbl{eqn.sum}
\sum_{i=0}^{I} a_i(n) \sum_{k}F(n+i, k) = 0.
\end{eqnarray}
In other words, we have:
\begin{equation}
\lbl{eq.haveS}
\sum_{i=0}^{I} a_i(n) S(n+i) = 0.
\end{equation}
So, Equation \eqref{eqn.celine} produces a recursion relation, which
is inhomogeneous if we are summing over an interval. How can we
find functions $a_i(n)$ that satisfy Equation~\eqref{eqn.celine}? The idea
is simple: divide Equation~\eqref{eqn.celine} by $F(n,k)$, and 
use \eqref{eq.rationalF} to convert the divided equation into a
{\em linear equation} over the field $\BQ(n,k)$, with unknowns
$a_i(n)$ for $i=0,\dots,I$. Clearing denominators, we get linear equation
over $\BQ(n)[k]$ with the same unknowns $a_i(n)$. 
Thus, the coefficients of every power of $k$ must vanish, and this gives a 
linear system of equations over $\BQ(n)$ with unknowns $a_i(n)$. If there are 
more unknowns than equations, one is guaranteed to find a nonzero solution. 
By a counting argument, one may see that if we choose $I$ high enough 
(this depends on the complexity of the term $F(n,k)$), then we have more 
equations than unknowns.

Although it can be numerically challenging to
find $a_i(n)$ that satisfy Equation~\eqref{eqn.celine}, it is routine
to check the equation once $a_i(n)$ are given. Indeed, one only need to divide 
the equation by $F(n,k)$, and then check that a function in $\BQ(n,k)$ is 
identically zero. The latter is computationally easy task in the field 
$\BQ(n,k)$.

This algorithm produces a recursion relation for $S(n)$. However, it is known 
that the algorithm does not always yield a recursion relation of the 
smallest order.

Applying Gosper's algorithm, Wilf and Zeilberger invented another algorithm, 
the {\em WZ algorithm}, also called {\em creative telescoping}. 
Instead of looking for $0$ on the right-hand side of 
Equation~\eqref{eqn.celine}, they instead looked for a function $G(n, k)$ 
such that

\begin{eqnarray} 
\lbl{eqn.wz}
\sum_{i=0}^{N} a_i(n) F(n+i, k) = G(n, k+1) - G(n, k).
\end{eqnarray}
Summing over $k$, and using telescoping cancellation of the terms
in the right hand side, we get a recursion relation for $S(n)$.
How to find the $a_i(n)$ and $G(n,k)$ that satisfy \eqref{eqn.wz}?
The idea is to look for a {\em rational function} $\Cert(n,k)$
(the so-called {\em certificate} of \eqref{eqn.wz}) such that

\begin{equation}
\lbl{eq.cert}
G(n, k) = \Cert(n, k)F(n, k).
\end{equation}
Dividing out \eqref{eqn.wz} by $F(n,k)$ and proceeding as before, 
one reduces this to a problem of linear algebra. As before, given $a_i(n)$ and
$\Cert(n,k)$, it is routine to check whether \eqref{eqn.wz} holds.

Now, let us rephrase the above equations using operators.
We define operators $E$, $E_k$, $n$ and $k$ that act on a function $F(n,k)$
by:

\begin{align}
\lbl{eq.Enk}
(E F)(n, k) &= F(n+1, k), &
(E_k F)(n, k) &= F(n, k+1), \\
(n F)(n,k) &=n F(n,k), & (kF)(n,k) &=kF(n,k).
\end{align}
The operators $E$ and $n$ (and also $E_k,k$) do not commute. Instead, we have:
$$
En=(n+1)E, \qquad E_k k=(k+1) E_k.
$$
On the other hand, $n,E$ commute with $k,E_k$. 
Then we can rewrite Equation~\eqref{eqn.wz} as

\begin{eqnarray} 
\lbl{eqn.operator}
\left(\sum_{i=0}^{I} a_i(n) E^i\right) F(n, k) = 
(E_k-1)G(n, k) = (E_k-1) \Cert(n, k) F(n, k).
\end{eqnarray}
Implementation of the algorithms are available in various platforms, 
such as,  Maple and Mathematica. See, for example, \cite{Z2} and \cite{PR2}.

Let us mention now how one deals with boundary terms.
In the applications below, one considers
not quite the unrestricted sums of Equation \eqref{eq.sumS}, but 
rather restricted ones of the form:

\begin{equation}
\lbl{eq.sumSp}
S'(n)=\sum_{k=0}^\infty F(n,k),
\end{equation}
where $F(n,k)$ is a proper hypergeometric term.
When we apply the Creating Telescoping summation to \eqref{eqn.wz}, 
we are left with some boundary terms $R(n) \in \BQ(n)$. In that case, 
Equation \eqref{eq.haveS} becomes:
$$
\left(\sum_{i=0}^{I} a_i(n) E^i\right) S'(n) = R(n).
$$
This is an inhomogeneous equation of order $I$ which can be converted
into a homogeneous recursion of order $I+1$ by following trick:
apply the operator 

$$
(E-1)\frac{1}{R(n)}
$$ 
on both sides of the recursion. We get

$$
\left(\frac{1}{R(n+1)}E - 
\frac{1}{R(n)}\right)\left(\sum_{i=0}^{I} a_i(n) E^i\right) S'(n) = 0,
$$
i.e.,

$$\left(\frac{a_{I}(n+1)}{R(n+1)}E^{I+1}
  + \sum_{i=1}^{I}\left(\frac{a_{i-1}(n+1)}{R(n+1)} 
- \frac{a_i(n)}{R(n)} \right) E^i 
- \frac{a_0(n)}{R(n)}\right) S'(n) = 0.
$$
In Quantum Topology we are using $q$-factorials rather than factorials.
The previous results translate without conceptual difficulty to the $q$-world,
although the computer implementation is slower.
A {\em term} $F(n, k)$ is called $q$-{\em hypergeometric} 
if 
$$
\frac{F(n+1, k)}{F(n, k)}, \frac{F(n, k+1)}{F(n, k)}
\in \BQ(q,q^n,q^k).
$$
Examples of $q$-hypergeometric terms are the quantum factorials of linear
forms in $n,k$, and ratios of products of quantum factorials
and $q$ raised to quadratic functions of $n$ and $k$.
The latter are called $q$-{\em proper hypergeometric}.

Sister Celine's algorithm and the WZ algorithm work equally well in the 
$q$-case. In either algorithms, we can replace the operators $E,n,E_k,k$
of \eqref{eq.Enk} by the operators $E,Q,E_k,Q_k$ defined by:

\begin{align}
\lbl{eq.Enkq}
(E F)(n, k) &= F(n+1, k), &
(E_k F)(n, k) &= F(n, k+1), \\
(Q F)(n,k) &=q^n F(n,k), & (Q_k F)(n,k) &=q^k F(n,k).
\end{align}
Observe that $E,Q$ (and also $E_K,Q_k$) $q$-commute, i.e., we have:

\begin{equation}
\lbl{eq.Encq}
EQ=qQE, \qquad E_k Q_k=qQ_kE_k.
\end{equation}
On the other hand, $E,Q$ commute with $E_k,Q_k$. With these modifications,
and with the replacement of the field $\BQ(n)$ by $\BQ(q,q^n)$,
the rest of the proofs still apply naturally.
The implementations of the $q$-case include \cite{PR3}, \cite{Ko} 
and \cite{Z2}.

\section{Multi-certificate creative telescoping and Theorem 
\ref{thm.1}}
\lbl{sec.hybrid}

\subsection{Multi-certificate Creative Telescoping}
\lbl{sub.hybrid}

In a nut-shell, the method of creative telescoping works as follows.
To find the recursion such that
$$
\sum^m_{i=0} a_i(n) \left(\sum_{k \ge 0} F(n+i, k)\right) = b(n),
$$
it suffices to find a rational function $\Cert(n,k) \in 
\BQ(q,q^n,q^k)$ such that $G(n, k):=\Cert(n,k)F(n,k)$ satisfies:
$$
\left(\sum^m_{i=0} a_i(n) E^i\right)F(n, k) 
= \sum^m_{i=0} a_i(n) F(n+i, k) = G(n, k+1) - G(n, k) = (E_k-1)G(n, k).
$$
If this can be done, we sum both sides for $0 \leq k < +\infty$, and
we obtain: 
$$
\sum^m_{i=0} a_i(n) \left(\sum_{k \ge 0} F(n+i, k)\right) = G(n, 0).
$$
For twist knots $K_p$, we have from Equation \eqref{eq.CJ}:

\begin{equation}
\lbl{eq.JJp}
J_p(n)=\sum_{k=0}^n c(n,k) \hatJ_p(k),
\end{equation}
where $c(n,k)$ is proper $q$-hypergeometric (given by \eqref{eq.cnk}),
and $\hatJ_p(n)$ satisfies a linear $q$-difference equation of degree $|p|$
from \cite{GS1}. 

Without loss of generality, suppose that $p>0$.
Suppose the minimal order recursion of $\hatJ_p(k)$ is:

\begin{equation}
\lbl{eq.recJhat}
\left(\sum^p_{i=0} r_i(k) E_k^i\right) \hatJ_p(k) = 0,
\end{equation}
with $r_p(k)=1$ and $r_i(k) \in \BQ(q,q^k)$ for $i=0,\dots,p$.
The idea is to look for $p$ certificates $\{C_0(n,k), \ldots C_{p-1}(n,k)\}$, 
such that 

\begin{equation}
\lbl{eq.manyc}
\left(\sum^m_{i=0} a_i(n) E^i\right)c(n, k)\hatJ_p(k) 
= (E_k-1)\left(\sum^{p-1}_{j=0}C_j(n,k)E_k^j\right)c(n, k)\hatJ_p(k).
\end{equation}

\subsection{A first reduction to linear algebra}
\lbl{sub.reduction1}

Our goal in this section, stated in Proposition \ref{prop.pcert} below,
is to translate the functional equation \eqref{eq.manyc}
into a system of linear equations with unknowns $a_i(n) \in \BQ(q,q^n)$
and $C_j(n,k) \in \BQ(q,q^n,q^k)$ for $i=0,\dots,m$ and $j=0,\dots,p-1$. 
Since $c(n,k)$ is proper $q$-hypergeometric, we have:
$$
\frac{Ec(n, k)}{c(n, k)} = s(n, k) \in \BQ(q,q^n,q^k), 
\quad \frac{E_kc(n, k)}{c(n, k)} = t(n, k) \BQ(q,q^n,q^k).
$$
Observe that
\begin{eqnarray*}
0 &=& \sum^p_{i=0} r_i(k)  \hatJ_p(k+i)		\\
  &=& \sum^p_{i=0} r_i(k)  
\frac{c(n, k+p)}{c(n, k+i)}\hatJ_p(k+i)c(n, k+i) \\  
&=& \sum^p_{i=0} r_i(k)  
\frac{c(n, k+p)}{c(n, k+i)}\hatJ_p(k+i)c(n, k+i) \\ 
&=& \left(\sum^p_{i=0} r_i(k)  
\frac{c(n, k+p)}{c(n, k+i)}E_k^i\right)\hatJ_p(k)c(n, k).
\end{eqnarray*}
So if we define 
\begin{eqnarray*}
R_i(n, k) &=& r_i(k) \frac{c(n, k+p)}{c(n, k+i)}  \\
&=& r_i(k) \prod^{p-i-1}_{j=0}t(n, k+i+j) \in \BQ(q,q^n,q^k),
\end{eqnarray*}
then we obtain that

\begin{equation}
\lbl{eq.Rnk}
\left(\sum^p_{i=0} R_i(n, k)E_k^i\right)\hatJ_p(k)c(n, k) = 0.
\end{equation}
Notice that since $r_p(k) = 1$, it follows that $R_p(n, k) = 1$ too.

\begin{proposition}
\lbl{prop.pcert}
Equation \eqref{eq.manyc} is equivalent to the following system
of linear equations:
\begin{equation}
\lbl{eq.prect1}
\sum^m_{i=0} a_i(n)\prod^{i-1}_{j=0}s(n+j, k) 
= -\sum^{p-1}_{j=0}C_{p-1}(n, k-j+1)R_{j}(n, k-j) 
\end{equation}
and
\begin{equation}
\lbl{eq.prect2}
C_{j-1}(n, k+1) = C_{j}(n, k) + C_{p-1}(n, k+1)R_{j}(n, k),
\qquad 1 \leq j \leq p-1,
\end{equation}
in the unknowns $a_i(n) \in \BQ(q,q^n)$ for $i=0,\dots,m$ and
$C_j(n,k) \in \BQ(q,q^n,q^k)$ for $j=0,\dots,p-1$. 
\end{proposition}

\begin{proof}
For convenience we define $C_{-1}(n, k) = 0$. Then using the
commutation relation 
$$
(E_k-1)C_j(n,k)=C_j(n,k+1)E_k - C_j(n, k)
$$
and
Equation \eqref{eq.Rnk} we obtain that:

\begin{eqnarray*}
\left(\sum^m_{i=0} a_i(n) E^i\right)c(n, k)\hatJ_p(k) 
&=& 
(E_k-1)\left(\sum^{p-1}_{j=0}C_j(n,k)E_k^i\right)c(n, k)\hatJ_p(k)	
\\ 
&=& 
\left(\sum^{p-1}_{j=0}C_j(n,k+1)E_k^{i+1} 
- \sum^{p-1}_{j=0}C_j(n,k)E_k^i\right)c(n, k)\hatJ_p(k)	
\\ 
&=& \left(\sum^{p}_{j=1}C_{j-1}(n,k+1)E_k^{i} 
- \sum^{p-1}_{j=0}C_j(n,k)E_k^i\right)c(n, k)\hatJ_p(k)	
\\ 
&=& \left(-C_{p-1}(n,k+1)\sum^{p-1}_{j=0}R_j(n, k)E_k^j \right. 
\\
& & \left.
+ \sum^{p-1}_{j=1}C_{j-1}(n,k+1)E_k^{i} 
- \sum^{p-1}_{j=0}C_j(n,k)E_k^i\right)c(n, k)\hatJ_p(k) 
\\ 
&=& \left(\sum^{p-1}_{j=0}
\left(-C_{p-1}(n,k+1)R_j(n, k)+C_{j-1}(n,k+1) \right.\right.
\\
& & \left.\left.
-C_j(n,k)\right)E_k^j\right)c(n, k)\hatJ_p(k).	
\end{eqnarray*}
So
\begin{eqnarray*}
& & \sum^m_{i=0} a_i(n) 
\left(\prod^{i-1}_{j=0}s(n+j, k)\right)c(n, k)\hatJ_p(k)
\\
&=& \left(\sum^m_{i=0} a_i(n) E^i\right)c(n, k)\hatJ_p(k)
\\ &=& \left(\sum^{p-1}_{j=0}
\left(-C_{p-1}(n,k+1)R_j(n, k)+C_{j-1}(n,k+1)-C_j(n,k)\right)E_k^j\right)
c(n, k)\hatJ_p(k).
\end{eqnarray*}
If we divide both sides by $c(n, k)$, which is hypergeometric, we obtain 
a new recursion on $\hatJ_p(k)$ of order $p-1$.
Since $J_p(k)$ satisfies a minimal order recursion of degree $p$, 
the last equality implies that
the coefficient of each $E_k^i$ is 0 for all $i$. Hence

$$
\left\{
\begin{array}{ll}-C_{p-1}(n,k+1)R_j(n, k)+C_{j-1}(n,k+1)-C_j(n,k) = 0	
& \mathrm{if\ } 1 \le j \le p-1, 
\\
\sum^m_{i=0} a_i(n)\frac{c(n+i, k)}{c(n, k)} 
= -C_{p-1}(n,k+1)R_0(n, k)-C_0(n,k)	& \mathrm{if\ } j=0.
\end{array}\right.
$$
The first equation implies \eqref{eq.prect2}. In particular,
$$
C_{0}(n, k+p-1) = C_{p-1}(n, k) 
+ \sum^{p-2}_{j=0}C_{p-1}(n, k+j+1)R_{p-j-1}(n, k+j).
$$
Therefore,
$$
\sum^m_{i=0} a_i(n)\frac{c(n+i, k)}{c(n, k)} 
= \sum^m_{i=0} a_i(n)\prod^{i-1}_{j=0}s(n+j, k) 
= -\sum^{p-1}_{j=0}C_{p-1}(n, k-j+1)R_{j}(n, k-j),
$$
which proves \eqref{eq.prect1} and concludes the proof of the proposition.
\end{proof}

\subsection{A second reduction to linear algebra}
\lbl{sub.reduction2}

Proposition \ref{prop.pcert} reduces the problem of finding $p$
certificates $C_j(n,k)$ to a problem of finding a single 
certificate $C_{p-1}(n,k)$. 
Since $C_{p-1}(n,k) \in \BQ(q,q^n,q^k)$ is a rational function, we can write
it in the form:

\begin{equation}
\lbl{eq.CND}
C_{p-1}(n, k) = \frac{N_p(n, k)}{D_p(n, k)},
\end{equation} 
where $N_p(n, k) = \sum^{r_d}_{i=0} d_i(n) q^{k i}$ and 
$D_p(n, k) = \sum^{r_e}_{i=0} e_i(n) q^{k i}$, 
in which $d_i(n)$ and $e_i(n)$ are in $\BQ(q, q^n)$.
By making the proper choices of $D_p(n, k)$ and the values of $m$ and $r_d$, 
we can clear denominators and 
convert Equation \eqref{eq.prect2} as a {\em linear equation}
in $\BQ(q,q^n)[q^k]$ with unknowns $a_i(n)$, $d_i(n)$, and $e_i(n)$.

Setting every coefficient of every power of $q^k$ to zero, we obtain a system
of linear equations in the unknowns $a_i(n)$, $d_i(n)$, and $e_i(n)$
and coefficients in the field $\BQ(q,q^n)$. A nontrivial solution is
guaranteed by Sister Celine's method for the case of $J_p(n)$. At any rate,
we can solve the system of equations using software like Maple or Mathematica.

Now comes the tricky part, and an educated
guess for the case of twist knots.
Since 
$$
\frac{c(n+i, k)}{c(n, k)} = \prod^{i-1}_{j=0}s(n+i, k),
$$  
and 
$R_{j}(n, k)$ are all polynomials in $\BQ(q, q^n)[q^k]$,
the most natural choice of $D_p(n, k)$ is the one such that $D_p(n, k-j+1)$ 
divides the polynomial $R_{j}(n, k-j)$ for all $j$.
Let $D_p(n, k)$, the denominator of the certificate $C_{p-1}(n, k)$, be
$$
\left\{
\begin{array}{lll}	
q^{pk} \prod^{p-1}_{i=1-p}(1 - q^{k-n-i})  & \mathrm{if\ } p > 0,
\\	
\prod^{2|p|}_{i=1}(1-q^{k-n+|p|-i})  & \mathrm{if\ } p < 0.
\end{array}
\right.
$$
While there is no guarantee that this will give
a minimal order inhomogeneous recursion relation for $J_p(n)$,
it does work for $-8 \leq p \leq 11$.
To show that the recursion relation for $(J_p)$ is of minimal order (and 
incidentally to check the AJ Conjecture of \cite{Ga}), we can set $q=1$.
An explicit computation shows that the $L$-degree of $A_p^{nh}(L,M,q)$
does not drop when we specialize to $q=1$ and moreover we have 
$$
A_p^{nh}(L,M,1)=A_p(L,M^2) F_p(M),
$$ 
where $A_p(L,M)$ is the $A$-polynomial of the twist knot $K_p$, and
$F_p(M) \in \BQ(M)$. 
$A_p(L,M^2)$ has been computed by Hoste-Shanahan in \cite{HS1}, and has been
shown to be irreducible in \cite{HS2}. It follows that $A_p^{nh}(E,Q,q)$
does not have any right factors and concludes the proof of Theorem \ref{thm.1}.
\qed

\section{Odds and ends}
\lbl{sub.odds}

\subsection{A generalization of Theorem \ref{thm.1}}
\lbl{sub.thm1}

In fact, the multi-certificate 
proof of Theorem \ref{thm.1} implies the following  result.

\begin{theorem}
\lbl{thm.11}
If $c(n,k)$ is proper $q$-hypergeometric term 
and $(\hatJ(n))$ is $q$-holonomic, and 
$$
J(n)=\sum_{k=0}^n c(n,k) \hatJ(k),
$$ 
then $J(n)$ is $q$-holonomic. A linear $q$-difference equation for $(J(n))$
can be constructed from a linear $q$-difference equation for $(\hatJ(n))$
and $c(n,k)$.
\end{theorem}
A software package that accompanies the proof Theorem \ref{thm.1} was
developed by the second author.

\begin{remark}
\lbl{rem.bigp}
Our proof of Theorem \ref{thm.1} reduces to solving a system of $2|p|$ linear 
equations
over the field $\BQ(q,q^n)$. When $-8 \leq p \leq 11$, this system
can be solved explicitly by symbolic software. Le's algorithm for
computing the non-commutative $A$-polynomial of a 2-bridge knot,
also requires a system of linear equations $(2|p|)!$ over the field 
$\BQ(q,q^n)$ in the case of twist knots; see \cite{Le}. However,
an implementation of Le's algorithm exceeded the capacity of our symbolic
software for $p=1$ and $p=-1$.
\end{remark}

\subsection{Is there a recursion of the non-commutative $A$-polynomial
with respect to the number of twists?}
\lbl{sub.rec}

Recall that $A_p(L,M)$ denotes the $A$-polynomial of the twist knot $K_p$.
In \cite{HS1}, Hoste-Shanahan use a trace identity in $\mathrm{SL}(2,\BC)$
in order to give a second order linear recursion relation for the sequence 
$(A_p)$.

There is supporting evidence that $A^{nh}_p(L,M,1)$ is annihilated by the
following operator

\begin{equation}
\lbl{eq.Apq}
\left\{
\begin{array}{lll}M^2(M-1)^4(M+1)^4(L+M)^4-(M-1)^2(M+1)^2(M^4-LM^4+2LM^3 & 
\\ \quad +L^2M^2+M^2+2LM^2+2LM+L^2-L)P+P^2 & \text{if } p > 0; \\
1-(M-1)^2(M+1)^2(M^4-LM^4+2LM^3+L^2M^2+M^2+2LM^2+2LM &  
\\ \quad +L^2-L)P+M^2(M-1)^4(M+1)^4(L+M)^4P^2 & \text{if } p < 0,
\end{array}
\right.
\end{equation}
where 

$$
P A^{nh}_p(L,M,1) = A^{nh}_{p+1}(L,M,1).
$$
Equation \eqref{eq.Apq} may be proven using the recursion on $\hatJ_p(k)$
and its simplification when $q=1$; see \cite[Thm.2]{GS1}. Unfortunately, 
there is equally strong evidence that the sequence $A^{nh}_p(E,Q,q)$
does not satisfy a linear recursion with respect to $p$.

\appendix

\section{A generating functions approach}
\lbl{sec.GF}

In this appendix we present an alternative approach to get a recursion relation
for $J_K(n)$ given Equation \eqref{eq.CJ} and a recursion relation for 
$\hatJ_K(n)$. This idea was communicated to us by D. Zeilberger, and may
be useful in its own right. We were not able to compute the non-commutative
$A$-polynomial for twist knots this way.

To explain the idea, let us recall first that a sequence $(a(n))$ of rational
numbers is holonomic iff the generating series 
$$
F(z)=\sum_{n=0}^\infty a(n) z^n
$$
is holonomic, i.e., it is annihilated by an element of the Weyl algebra
$\BQ\la z, d/dz \ra$; see \cite{Z1}. The $q$-analogue of this is the
following. Consider a sequence $(a(n))$ with $a(n) \in \BQ(q)$, and the
generating series

\begin{equation}
\lbl{eq.Fzq}
F(z,q)=\sum_{n=0}^\infty a(n) z^n \in \BQ(q)[[z]]
\end{equation}
There are two operators $Q$ and $Z$ that act on the elements $F(z,q)$ of
$\BQ(q)[[z]]$ by:

\begin{equation*}
(QF)(z,q)=F(qz,q), \qquad (ZF)(z,q)=zF(z,q)
\end{equation*}
It is easy to see that $QZ=qZQ$, and that $(a(n))$ is $q$-holonomic iff
the generating series $F(z,q)$ is $q$-holonomic.

Now, let us consider two sequences $(J(n))$ and $(\hatJ(n))$ of rational
functions that are related by:

\begin{equation}
\lbl{eq.CJ2}
J(n) = \sum^n_{k=0} c(n, k) \hatJ(k)
\end{equation}
where the kernel $c(n,k)$ is given by \eqref{eq.cnk}.
$c(n, k)$ can be slightly simplified into 
$q^{-nk}\frac{\qfac{n+k}}{\qfac{n-k-1}(1-q^n)}$. We will absorb the
factor $\frac{1}{1-q^n}$ in the colored Jones function and define

\begin{eqnarray*}
\gamma(n,k) &:=& q^{-nk}\frac{\qfac{n+k}}{\qfac{n-k-1}} \\
&=&
c(n,k) (1-q^n) \\
\checkJ(n) &:=& \sum_{k=0}^n \gamma(n,k) \hatJ(k) \\
&=& \hatJ(n) (1-q^n).
\end{eqnarray*}

\begin{proposition}
\lbl{prop.Hkz}
Let  
\begin{equation}
\lbl{eq.Hkz}
H(k,z) = \sum^{\infty}_{i=0} \gamma(k+i,k) z^{i},
\end{equation}
then we have:
\begin{equation}
\lbl{eq.Hkz2}
H(k,z) = (-1)^k q^{\frac{-k(k+1)}{2}}
\frac{\qfac{2k+1}}{z^k\qrf{\frac{q}{z}}{k}\qrf{z}{k+2}}.
\end{equation}
\end{proposition}

\begin{proof}
We will use the idea of the WZ-algorithm to the hypergeometric summand:

\begin{equation} 
\lbl{eqn.gf.summand}
H_1(k, i, z) = \gamma(k+i,k) z^i.
\end{equation}
We claim that:

\begin{equation}
\lbl{eq.H1a}
(z-q^{k+1})H_1(k+1, i, z) + 
\frac{(1-q^{2k+2})(1-q^{2k+3})}{q^{k+1}(1-zq^{k+2})}H_1(k, i, z) 
= G_1(k, i, z) - G_1(k, i-1, z),
\end{equation}
where

\begin{equation}
\lbl{eq.H1b}
G_1(k, i, z) = 
\frac{-z(1-q^{2k+i+1})
(zq^{k+2}-1-zq^{3k+i+4}+q^{4k+i+5}}{q^{2k+i+1}(1-zq^{k+2})}H_1(k, i, z).
\end{equation}
Equation \eqref{eq.H1a} can be verified by 
dividing both sides by $H_1(k, i, z)$, and then it reduces to an
identity in the field $\BQ(z,q,q^k)$ which can be readily checked.
Now summing both sides of Equation \eqref{eqn.gf.summand} over $k$, 
and we get the desired result.
\end{proof}

Consider the generating function of $\checkJ(n)$: 

\begin{equation}
\lbl{eq.Fpz}
F(z,q) = \sum^{\infty}_{n=0} \checkJ(n) z^n.
\end{equation}

\begin{proposition}
\lbl{prop.Fzq}
We have:
\begin{equation}
\lbl{eq.FJzq}
F(z,q) = \sum^{\infty}_{k=0} (-1)^k q^{\frac{-k(k+1)}{2}}
\frac{\qfac{2k+1}}{\qrf{\frac{q}{z}}{k}\qrf{z}{k+2}} \hatJ_p(k).
\end{equation}
\end{proposition}

\begin{proof}
We will interchange the order of summation and use Proposition 
\ref{prop.Hkz}. We get: 

\begin{eqnarray*}
F(z,q)	
&=& \sum^{\infty}_{n=0} \checkJ(n) z^n	\\
&=& \sum^{\infty}_{n=0} \left(\sum^n_{k=0} \gamma(n,k) \right) \hatJ_p(k) z^n
\\
&=& \sum^{\infty}_{k=0} \hatJ_p(k)  \left(\sum^{\infty}_{n=k} \gamma(n,k) 
z^{n} \right),	
\\
&=& \sum^{\infty}_{k=0} \hatJ_p(k)  z^k \left(\sum^{\infty}_{i=0} 
\gamma(k+i,k) z^{i} \right),	
\\
&=& \sum^{\infty}_{k=0} H(k,z) \hatJ_p(k)  z^k	\\
&=& \sum^{\infty}_{k=0} (-1)^k 
q^{\frac{-k(k+1)}{2}}
\frac{\qfac{2k+1}}{\qrf{\frac{q}{z}}{k}\qrf{z}{k+2}} \hatJ(k).
\end{eqnarray*}
\end{proof}
One can use a $q$-difference equation for $\hatJ(n)$ and
Proposition \eqref{eq.Fzq} to get a $q$-difference equation for $F(z,q)$.
This will be explored in another publication.

\section{The non-commutative $A$-polynomial for $p=\pm 3$}
\lbl{sec.p-33}

In this section we give the inhomogeneous non-commutative $A$-polynomial
for $p=\pm 3$. The reader may compare the size of the output with Theorem
\ref{thm.0}.


{\small  
\begin{eqnarray*}
A^{nh}_3(E,Q,q)	
&=& q^{11n+20}(-1+q^n)(q^{2n+8}-1)(q^{2n+7}-1)(q^{2n+6}-1)(q^{2n+9}-1)\\
& & -q^{4n+14}(q^{n+1}-1)(q^{2n+2}-1)(q^{2n+7}-1)(q^{2n+8}-1)(q^{2n+9}-1)\\
& & (-1-q^{3+2n}+q^{8+4n}+q^{1+n}+q^{4+3n}+q^{3n+3}+q^{9+4n}-q^{5+4n}-q^{7+4n}\\
& & -q^{10+5n}-q^{8+6n}+q^{4n+3}+q^{1+2n}-q^{6+4n}-q^{9+5n}-q^{2n+2}-q^{10+6n}\\
& & -q^{7+3n}-q^{5n+4}+q^{11+7n}-2q^{9+6n}+q^{10+7n}-q^{5n+5}+q^{6+2n}+q^{8+5n}\\
& & -q^{2+3n})E+q^{2n+9}(q^{n+2}-1)(q^{2n+1}-1)(q^{2n+4}-1)(q^{2n+8}-1)\\
& & (q^{2n+9}-1)(-1-2q^{13+6n}-q+q^{3+2n}-q^{8+4n}-q^{17+8n}+2q^{16+7n}-2q^{4+3n}\\
& & +2q^{12+4n}-2q^{3n+9}-q^{18+7n}-2q^{18+8n}-4q^{9+4n}-q^{14+7n}-2q^{10+4n}\\
& & -q^{19+8n}+2q^{2+n}+2q^{7+4n}+q^{10+5n}-q^{11+3n}+q^{15+6n}+q^{8+3n}\\
& & +4q^{5n+12}+q^{17+6n}+q^{8+2n}-2q^{5+2n}-2q^{10+3n}+q^{7+5n}+q^{17+7n}\\
& & -q^{14+4n}+2q^{6+4n}-q^{9+5n}-q^{13+7n}+q^{2n+2}-2q^{14+6n}+q^{10+6n}\\
& & +q^{4+n}+5q^{11+5n}+3q^{7+3n}+2q^{15+7n}+2q^{3+n}-q^{8+5n}-2q^{2n+4}+2q^{11+4n}\\
& & +2q^{7+2n}+2q^{11+6n}+q^{20+9n}-q^{15+5n}-2q^{5+3n}+3q^{6n+16}-2q^{14+5n}\\
& & +q^{6+3n}+q^{12+6n})E^2-q^3(q^{n+3}-1)(q^{2n+1}-1)(q^{2n+2}-1)(q^{2n+6}-1)\\
& & (q^{2n+9}-1)(-1+q^{3+2n}+q^{8+4n}-4q^{12+4n}+2q^{3n+9}-q^{18+7n}+q^{9+4n}\\
& & -2q^{8n+22}-q^{10+4n}-q^{24+8n}-2q^{7n+22}-q^{18+6n}+q^{2+n}-q^{7+4n}\\
& & -2q^{17+5n}-3q^{11+3n}+q^{6n+21}+2q^{15+6n}+2q^{8+3n}-2q^{5n+12}+2q^{20+7n}\\
& & -3q^{17+6n}+q^{8+2n}-2q^{5+2n}-q^{10+3n}-q^{17+7n}+2q^{14+4n}+q^{4n+15}\\
& & +q^{9n+25}+2q^{14+6n}+q^{4+n}+2q^{19+7n}-2q^{11+5n}-q^{7+3n}-2q^{16+5n}\\
& & +2q^{6n+20}-q^{23+7n}-2q^{6+2n}+2q^{3+n}+q^{2n+4}+q^{5n+13}-5q^{11+4n}\\
& & -q^{7+2n}+2q^{6n+19}+2q^{15+5n}-q^{5+3n}-q^{3n+12}+q^{19+5n}-q^{6n+16}\\
& & +4q^{14+5n}-2q^{23+8n}+q^{9n+26}-2q^{6+3n})E^3+q(q^{n+4}-1)(q^{2n+1}-1)\\
& & (q^{2n+2}-1)(q^{2n+3}-1)(q^{2n+8}-1)(-1-q+q^{12+4n}+q^{4n+17}+q^{2n+9}+q^{17+5n}\\
& & +q^{11+3n}-q^{6n+21}-q^{8+3n}-q^{8+2n}+q^{18+5n}+q^{5+2n}+q^{10+3n}-q^{14+4n}\\
& & -q^{13+3n}-q^{5n+21}+2q^{4+n}+q^{25+7n}-q^{16+5n}+q^{3+n}+q^{2n+4}-q^{14+3n}\\
& & -q^{13+4n}+q^{n+5}+q^{3n+12}+q^{2n+10})E^4+(q^{n+5}-1)(q^{2n+1}-1)(q^{2n+2}-1)\\
& & (q^{2n+3}-1)(q^{2n+4}-1)E^5,
\end{eqnarray*}
\begin{eqnarray*}
B_3(q^n,q)
&=& q^{3n+12}(q^{n+1}+1)(q^{n+2}+1)(q^{n+3}+1)(q^{n+4}+1)(q^{2n+1}-1)(q^{2n+3}-1)\\
& & (q^{2n+5}-1)(q^{2n+7}-1)(q^{2n+9}-1),
\end{eqnarray*}
\begin{eqnarray*}
A^{nh}_{-3}(E,Q,q)
&=& q^{6n+18}(-1+q^n)(q^{2n+7}-1)(q^{2n+8}-1)(q^{2n+9}-1)(q^{2n+10}-1)(q^{2n+11}-1)\\
& & -q^{4n+14}(q^{n+1}-1)(q^{2n+2}-1)(q^{2n+8}-1)(q^{2n+9}-1)(q^{2n+10}-1)(q^{2n+11}-1)\\
& & (1+q+q^{5n+12}+q^{4+2n}-q^{2n+7}-q^{13+7n}+q^{8+4n}-q^{2n+2}+q^2-q^{5+n}+q^{8n+14}\\
& & -q^{8+3n}-q^{2n+9}-q^{4+n}-q^{2+n}-q^{7+6n}-q^{8+2n}-q^{11+5n}-q^{1+2n}+q^{4+3n}\\
& & -q^{7+3n}-q^{10+5n}+q^{11+6n}-q^{11+4n}+q^{9+4n}-q^{5+4n}+q^{6+5n}-2q^{3+n}+q^{12+6n}\\
& & -q^{13+6n}+q^{3n+3}-q^{5+3n}-q^{4n+4}+q^{6+4n}+q^{10+3n})E+q^{2n+9}(q^{n+2}-1)\\
& & (q^{2n+1}-1)(q^{2n+4}-1)(q^{2n+9}-1)(q^{2n+10}-1)(q^{2n+11}-1)(1-q^{9n+24}+3q^{10+4n}\\
& & +q+q^{16+7n}-3q^{8+4n}-q^{4+2n}-3q^{4+n}-2q^{6n+19}+4q^{5+3n}+3q^{15+7n}-3q^{3n+9}\\
& & +q^{4+3n}-q^{11+6n}+q^{19+8n}+q^2-q^{18+6n}+q^{24+10n}+q^{4n+17}-2q^{12+6n}+2q^{14+7n}\\
& & +q^{25+10n}+q^{10+6n}+3q^{17+5n}-q^{2n+9}-2q^{9n+23}+3q^{6n+16}-2q^{3+2n}-4q^{14+5n}\\
& & +q^{18+5n}-q^{2n+10}-q^{2+n}+q^{17+6n}+2q^{11+5n}+2q^{8n+21}-q^{23+8n}-2q^{8+3n}\\
& & -q^{2+2n}-2q^{17+7n}+3q^{10+5n}+3q^{3n+12}+4q^{15+6n}+2q^{5+2n}-q^{21+9n}+2q^{9+5n}\\
& & -q^{7+5n}+2q^{7+2n}+q^{7n+22}-2q^{5n+12}+2q^{8n+20}-3q^{18+7n}-q^{8n+16}+3q^{6+3n}\\
& & +q^{5+4n}-q^{6+n}+q^{13+4n}-2q^{17+8n}+6q^{11+4n}-2q^{5+n}+q^{9+6n}+2q^{13+3n}\\
& & -6q^{5n+13}-3q^{3+n}+4q^{12+4n}-2q^{19+7n}-q^{6n+20}-3q^{7+4n}+q^{14+3n}+2q^{6+2n}\\
& & +2q^{16+5n}-q^{18+8n}+3q^{11+3n}-3q^{13+6n}+2q^{7n+21}+q^{20+7n}-2q^{22+9n}\\
& & +q^{13+7n})E^2-q^3(q^{n+3}-1)(q^{2n+1}-1)(q^{2n+2}-1)(q^{2n+6}-1)(q^{2n+10}-1)\\
& & (q^{2n+11}-1)(1+3q^{9n+24}-4q^{10+4n}-2q^{16+7n}-q^{8+4n}-q^{4+2n}-q^{4+n}+7q^{6n+19}\\
& & +q^{5+3n}-q^{15+7n}-3q^{3n+9}+q^{19+8n}+4q^{18+6n}-q^{28+10n}+q^{36+12n}-q^{4n+17}\\
& & +2q^{14+4n}+7q^{25+8n}-3q^{10+3n}-3q^{17+5n}-4q^{6n+21}+q^{9n+23}-8q^{6n+16}-q^{3+2n}\\
& & -q^{14+5n}-q^{2+n}-q^{17+6n}+q^{8+2n}+q^{11+5n}-3q^{8n+21}+q^{23+8n}-3q^{23+6n}\\
& & +q^{10n+32}-q^{27+8n}-q^{18+4n}-2q^{14+6n}+3q^{31+10n}-q^{7n+28}+q^{17+7n}-q^{34+11n}\\
& & -2q^{10+5n}+q^{3n+12}-2q^{33+11n}-6q^{15+6n}+q^{5+2n}-q^{9+5n}+q^{9n+32}+2q^{26+8n}\\
& & -3q^{9+4n}-3q^{27+9n}-4q^{8n+22}-3q^{8n+28}+3q^{7+2n}-q^{11n+32}-8q^{7n+22}-2q^{24+6n}\\
& & +4q^{5n+12}-q^{8n+20}+4q^{18+7n}+2q^{9n+31}-3q^{9n+28}-q^{8n+30}+3q^{6+3n}-q^{22+5n}\\
& & -q^{8n+29}-q^{4n+15}+7q^{13+4n}+q^{11+4n}-q^{5+n}+2q^{13+3n}+5q^{5n+13}+6q^{24+8n}\\
& & -2q^{3+n}+6q^{12+4n}+5q^{19+7n}+3q^{10n+30}-6q^{22+6n}+q^{7+4n}+3q^{9n+25}+q^{14+3n}\\
& & +3q^{6+2n}+2q^{19+5n}-8q^{16+5n}-q^{10n+27}-6q^{15+5n}+q^{13+6n}+q^{9n+30}-3q^{4n+16}\\
& & +3q^{7+3n}-6q^{7n+21}-q^{20+7n}+q^{29+10n}-q^{35+11n}-3q^{23+7n}+2q^{25+7n})E^3\\
& & +q(q^{n+4}-1)(q^{2n+1}-1)(q^{2n+2}-1)(q^{2n+3}-1)(q^{2n+8}-1)(q^{2n+11}-1)(1+4q^{24+6n}\\
& & +q^{10+4n}+q-q^{33+8n}+3q^{30+7n}-q^{4+2n}-q^{26+8n}+q^{36+10n}-2q^{4+n}+2q^{31+7n}\\
& & +2q^{15+5n}-3q^{3n+12}+3q^{3n+9}-4q^{20+5n}+q^{37+10n}+4q^{23+7n}-3q^{7n+27}+q^{25+6n}\\
& & +q^{32+7n}-2q^{18+5n}+q^{4n+17}+3q^{23+5n}+q^{10+3n}-2q^{26+7n}+q^{38+10n}-q^{2n+11}\\
& & +3q^{22+6n}+q^{17+6n}+q^{24+5n}+3q^{24+7n}+2q^{8+2n}-q^{6+n}+q^{29+6n}-2q^{27+8n}\\
& & -q^{18+4n}+2q^{15+3n}-q^{34+8n}-2q^{13+3n}-3q^{9n+34}-3q^{6n+20}-2q^{5+2n}-q^{8n+28}\\
& & -q^{9n+32}+q^{16+3n}+q^{9+4n}+q^{7+2n}+q^{7n+22}-6q^{19+5n}-q^{36+9n}+2q^{2n+9}\\
& & +2q^{8n+30}-3q^{9n+33}+6q^{23+6n}+3q^{16+5n}+2q^{22+5n}+2q^{8n+29}+4q^{4n+15}-3q^{13+4n}\\
& & -2q^{4n+19}-q^{11+4n}-2q^{5+n}-q^{5n+13}-q^{3+n}+3q^{7n+29}-2q^{12+4n}+3q^{4n+16}\\
& & +2q^{17+5n}+q^{14+3n}-q^{20+4n}-q^{6+2n}-2q^{11+3n}-3q^{6n+19}+2q^{8n+31}+q^{7+3n}\\
& & +2q^{8+3n}-2q^{35+9n})E^4-(q^{n+5}-1)(q^{2n+1}-1)(q^{2n+2}-1)(q^{2n+3}-1)(q^{2n+4}-1)\\
& & (q^{2n+10}-1)(1-q^{15+3n}+q^{18+4n}-q^{31+6n}+q^{4n+19}-q^{23+6n}+q^{16+3n}-q^{14+3n}\\
& & +q^{2n+9}+q^{2n+10}+q^{34+8n}+q^{20+5n}-q^{24+6n}+q^{26+5n}-q^{30+7n}-q^{32+7n}+q^{36+8n}\\
& & +q^{26+6n}-q^{5+2n}+q^{19+5n}-q^{29+6n}-q^{6n+30}-q^{21+4n}-q^{7n+33}-q^{23+5n}-q^{5+n}\\
& & -q^{24+5n}+q^{35+8n}-q^{11+2n}-q^{5n+21}-2q^{31+7n}+q^{4n+16}-q^{4n+15}+q^{10+3n}\\
& & -q^{14+4n})E^5+q^{6n+30}(q^{n+6}-1)(q^{2n+1}-1)(q^{2n+2}-1)(q^{2n+3}-1)(q^{2n+4}-1)\\
& & (q^{2n+5}-1)E^6,
\end{eqnarray*}
\begin{eqnarray*}
B_{-3}(q^n,q)&=& 
q^{3n+15}(q^{n+1}+1)(q^{n+2}+1)(q^{n+3}+1)(q^{n+4}+1)(q^{n+5}+1)(q^{2n+1}-1)(q^{2n+3}-1)\\
& & (q^{2n+5}-1)(q^{2n+7}-1)(q^{2n+9}-1)(q^{2n+11}-1).
\end{eqnarray*}
}
For a computer data of the non-commutative $A$-polynomial of twist
knots, see \cite{GS2}.

\section{The $A$-polynomial for $p=-3,\dots,3$}
\lbl{sec.Acomm}

For comparison, we give a formula of the $A$-polynomial $A_p(L,M)$
of the twist knot $K_p$, taken from \cite{HS1}.


\begin{center}
\begin{tabular}{|c|l|} \hline
$p$ & $A_p(L,M)$ 
\\ \hline
$1$ & $L + M^6$ 
\\ \hline
$-1$ & $-L + L M^2 + M^4 + 2 L M^4 + L^2 M^4 + L M^6 - L M^8$ 
\\ \hline
$2$ & $-L^2 + L^3 + 2 L^2 M^2 + L M^4 + 2 L^2 M^4 - L M^6 - L^2 M^8 + 
 2 L M^{10} + L^2 M^{10} + 2 L M^{12} + M^{14} - L M^{14} $ 
\\ \hline
$-2$ & $ L^2 - L^3 - 3 L^2 M^2 + L^3 M^2 - 2 L M^4 - L^2 M^4 + 3 L M^6 + 
 3 L^2 M^6 + M^8 + 3 L M^8 + 6 L^2 M^8 + 3 L^3 M^8 
$ 
\\ \hline
& 
$
+ L^4 M^8 + 
 3 L^2 M^{10} + 3 L^3 M^{10} - L^2 M^{12} - 2 L^3 M^{12} + L M^{14} - 
 3 L^2 M^{14} - L M^{16} + L^2 M^{16}$ 
\\ \hline
$3$ & $L^3 - 2 L^4 + L^5 - 4 L^3 M^2 + 4 L^4 M^2 - 2 L^2 M^4 + 2 L^3 M^4 + 
 3 L^4 M^4 + 5 L^2 M^6 + 5 L^3 M^6 + L M^8 
$ 
\\ \hline
& 
$
+ L^2 M^8 + 6 L^3 M^8 - 
 L M^{10} - 4 L^2 M^{10} - 4 L^3 M^{12} - L^4 M^{12} + 6 L^2 M^{14} 
+ L^3 M^{14} +
  L^4 M^{14} + 5 L^2 M^{16} 
$ 
\\ \hline
& 
$
+ 5 L^3 M^{16} + 3 L M^{18} + 2 L^2 M^{18} - 
 2 L^3 M^{18} + 4 L M^{20} - 4 L^2 M^{20} + M^{22} - 2 L M^{22} + L^2 M^{22} $ 
\\ \hline
$-3$ & $-L^3 + 2 L^4 - L^5 + 5 L^3 M^2 - 6 L^4 M^2 + L^5 M^2 + 3 L^2 M^4 - 
 6 L^3 M^4 - 10 L^2 M^6 - 5 L^3 M^6 + 4 L^4 M^6 
$ 
\\ \hline
& 
$
- 3 L M^8 - 
 3 L^3 M^8 + 5 L M^{10} 
+ 12 L^2 M^{10} + 10 L^3 M^{10} + M^{12} + 4 L M^{12} + 
 10 L^2 M^{12} + 20 L^3 M^{12} 
$ 
\\ \hline
& 
$
+ 10 L^4 M^{12} 
+ 4 L^5 M^{12} + L^6 M^{12} + 
 10 L^3 M^{14} + 12 L^4 M^{14} + 5 L^5 M^{14} - 3 L^3 M^{16} 
- 3 L^5 M^{16} + 
 4 L^2 M^{18} 
$ 
\\ \hline
& 
$
- 5 L^3 M^{18} - 10 L^4 M^{18} 
- 6 L^3 M^{20} + 3 L^4 M^{20} + 
 L M^{22} - 6 L^2 M^{22} + 5 L^3 M^{22} 
- L M^{24} + 2 L^2 M^{24} - L^3 M^{24} $ 
\\ \hline
\end{tabular}
\end{center}

\ifx\undefined\bysame
        \newcommand{\bysame}{\leavevmode\hbox
to3em{\hrulefill}\,}
\fi

\end{document}